\documentclass[11pt,reqno,draft]{amsart}
\usepackage{latexsym,amssymb}

\newtheorem{defi}{Definition}[section]

\newtheorem{cor}[defi]{Corollary}

\newtheorem{rem}[defi]{Remark}

\newtheorem{theorem}[defi]{Theorem}

\newcommand{\defeq}{\mathrel{\mathrm{\raise0.1ex\hbox{:}\hbox{=}\strut}}}

\def\dd{\:\mathrm{d}}

\def\R{\mathbb R}

\def\E{\mathbb E}

\def\F{\mathcal{ F}}

\def\ve{\varepsilon}
\def\C{\mathcal C}
\def\P{\mathbb P}
\begin{document}

\numberwithin{equation}{section}

\title[Harnack Inequality for Memory  SDEs]{Harnack Inequality for Functional SDE{\tiny s} with Bounded Memory}
\author[A. Es--Sarhir]{Abdelhadi Es--Sarhir}
\author[M-K. von Renesse]{Max-K. von Renesse}
\author[M. Scheutzow]{Michael Scheutzow}

\address{Technische Universit\"at Berlin, Fakult\"at II, Institut f\"ur Mathematik, Sekr. Ma 7-4\newline Stra{\ss}e des 17. Juni 136, D-10623 Berlin, Germany}

\email{[essarhir,mrenesse,ms]@math.tu-berlin.de}

\def\krylovrozovskii{Krylov-Rozovski{\u\i}}
\def\subvstern{_{V^*}}

\keywords{Harnack Inequality, Coupling, Strong Feller Property}


\begin{abstract}
We use a  coupling method for functional stochastic differential equations with bounded memory to establish an analogue of Wang's dimension-free Harnack inequality \cite{MR1481127}. The strong Feller property  for the corresponding segment process is also obtained.  \end{abstract}

\maketitle
\section{Introduction and Statement of Results}

\noindent Harnack inequalities are known to hold for a wide range of Markov processes such as diffusions or symmetric jump processes on manifolds, graphs, fractals or even more general metric measure spaces, where they  are a fundamental tool for the analysis of the corresponding transition semigroups, cf. \cite{MR2291922}. In most cases, a Harnack inequality is established under appropriate ellipticity conditions by harmonic analysis arguments. Such  arguments are typically not applicable in  non-Markovian or infinite dimensional set-ups. However, as shown in \cite{MR2215664} for finite dimensional diffusions, the dimension-free Harnack inequality of Wang \cite{MR1481127} may be proved by a purely probabilistic coupling technique which was recently adapted to the infinite dimensional case of monotone SPDEs \cite{MR2330974,liu-2008,DaPrato2009992}. \smallskip

In this note we show that the coupling method  works well also in the non-Markovian case of stochastic functional equations with additive noise and  Lipschitz drift with bounded delay. Moreover, the strong Feller property is obtained for the corresponding infinite dimensional segment process. The main additional difficulty compared to the diffusive case \cite{MR2215664} is the necessity to couple together two solutions including their pasts  at a given time. It turns out that this can be done simply by driving the second process with the drift induced from the segment process of the first.
\smallskip

For a precise statement of our results, fix  $r>0$ and  let $\C$ denote the space of continuous $\R^d$-valued functions
on $[-r,0]$ endowed with the sup-norm $\|\cdot \|$. For a function or a
process $X$ defined on $[t-r,t]$ we write $X_t(s):=X(t+s)$, $s \in
[-r,0]$.  Consider the stochastic functional differential equation
\begin{equation}
\label{sfde}
\left\{
\begin{aligned}
\dd X(t) & = V(X_t)\,\dd t + \dd W(t),\\
X_0& = \varphi,\\
 \end{aligned}
\right.
\end{equation}
where $W$ is an $\R^d$-valued Brownian motion  defined on a complete probability space $(\Omega,\F,\P)$ endowed with the augmented Brownian
filtration $\F^W_t = \sigma\left(W(u), {0\leq u \leq t} \right)\vee \mathcal N\subset \F$, where $\mathcal N$ denotes the null-sets in
$\F$, $\varphi$ is an $(\F^W_t)$-independent  $\C$-valued random variable and  $V:\C \to \R^d$ is  a measurable map.
  \smallskip

Below we assume that $V$ admits a decomposition
\begin{equation}
 V(x) = v(x(0))  + Z(x),  \label{Vdecomp}
\end{equation}
 where $v\in C(\R^d;\R^d)$ is a dissipative vector field on $\R^d$, i.e.
\begin{equation} \langle v(a)-v(b), a-b\rangle \leq 0  \quad \forall \, a,b \in \R^d\label{weakdissip} \end{equation}
and $Z$ is globally Lipschitz on $\C$, i.e.\ for some $L>0$,
\begin{equation}\label{lipcon}
 \begin{minipage}{12cm}
\[| Z(x)-Z(y)| \leq L\, \|x-y\| \quad \forall x,y \in \C.\]
 \end{minipage}
\end{equation}

Global existence and uniqueness for  \eqref{sfde}  hold under much weaker conditions, c.f.\ \cite{Renesse:arXiv0812.1726}. In particular, the corresponding segment process $\{X_t^\varphi \in \C\,|\, t\geq 0, \varphi \in \C, t\mapsto X^\varphi(t) \mbox{ solves } \eqref{sfde}\}$ induces a Markov semigroup $(P_t)$ on $\C$ via $\varphi \mapsto  P_tf(\varphi)  = \E (f(X_t^\varphi))$,   for bounded measurable $f : \C\to \R$.\smallskip

 Now our main result is the following version of Wang's dimension free Harnack inequality for the semigroup $(P_t)$.

\begin{theorem} \label{mainthm} Assume that  $V=v+Z$ as in \eqref{Vdecomp} with dissipative $v$ and $\mbox{Lip}_\C(Z)\leq L$, then for any $p>1, T>r$ and any bounded measurable $f: \C \rightarrow \R$,
\begin{equation}
(P_T f(y) )^p \leq P_T(f^p)(x) \exp\Bigl( \frac{ p}{p-1}
\rho_T^2(x,y)\Bigr) \quad \forall \, x, y \in \C, \label{harnack}
\end{equation}
where
\[\rho_T^2(x,y) = \inf_{s \in ]r,   T]} \left\{\frac  {|x(0)-y(0)|^2} {s-r} +  s\, L^2\|x-y\|^2\right\}.\]
\end{theorem}

\begin{rem} {\normalfont Elementary computation yields
 \begin{equation*}
\rho^2_T(x,y)  =  \left\{
\begin{array}{ll}
  \frac  {|x(0)-y(0)|^2} {T-r} +  T L^2\|x-y\|^2 & \mbox{for } T\leq r + \frac{|x(0)-y(0)|}{L\, \|x-y\|} \\
  2 L |x(0)-y(0)| \cdot \| x -y\| + r \, L^2 \|x-y\|^2 & \mbox{for } T \geq   r + \frac{|x(0)-y(0)|}{L\, \|x-y\|}.
\end{array}
\right.
\end{equation*}}
\end{rem}

Moreover, $(X_t)$ exhibits the following strong Feller property on the infinite dimensional state space $\C$.
Since the driving noise for $(X_t)$ is only $d$\,-\,dimensional,
this is a  noticeable result.

\begin{cor} \label{stronfell} Under  \eqref{Vdecomp}-\eqref{lipcon} the segment process $(X^x_t)$ is eventually strong Feller, i.e. let  $f: \C \to \R$ be bounded measurable, then for $t > r$ the  map $  x \mapsto P_tf(x)\in \R$ is continuous on $\C$.
\end{cor}

%
\begin{rem} {\normalfont
Our proofs below can easily be modified to include the  case of
random $V$ and random, strictly elliptic diffusion coefficient
$\sigma=\sigma(t)\in \R^{d\times d}$ in front of $\dd W(t)$ in
\eqref{sfde}. However, $(X_t)$ can generally not be expected to be
strong Feller in case $\sigma$ depends on the segment $X_t$, i.e.\
$\sigma=\sigma(X_t)$. If, for example, $d=1$ and the diffusion part
in \eqref{sfde} is of the form $g(X(t-1))$ or $g(\int_{t-r}^t
X(s)\dd s)$ with smooth, strictly increasing and positive $g$, then
the transition probabilities $P_t(x,\dd y)$ and $P_t(z, \dd y)$ are
mutually singular for all $t>0$ whenever $x\ne z$, since the initial
condition can perfectly recovered from $X_t$,  using the law of the
iterated logarithm c.f.\ \cite{MR2147280}.}
\end{rem}

Another  straightforward consequence of Theorem \ref{mainthm} is the following smoothing property of $(P_t)$.  For more on this we refer to \cite{MR2330974, DaPrato2009992}.
\begin{cor} Assume that $(X_t)$ admits  an invariant measure  $\mu \in \mathcal M(\C)$  such that
\begin{equation}\int_\C e^{\lambda \|x\|^2 } \mu(\dd x) < \infty   \mbox{ for some } \lambda > 4(2L+rL^2), \label{intcond}  \end{equation}
then for $t > r+ L^{-1}$,  $P_t$ is $\mu$-hyperbounded  i.e.\ $P_t$
is a bounded operator from $L^2(\C,\mu)$ to $L^4(\C,\mu)$.
\end{cor}

\noindent In the following we give an example when \eqref{intcond}
holds. For $z\in\R^d$ we set $v(z)=-\lambda_0 z$ for some
$\lambda_0>4(2L+rL^2)$. Furthermore, assume that $\sup\limits_{x\in
\R^d}\|Z(x)\|\leq M$ for some $M\geq 0$. Clearly in this case the
solution $(X(t))_{t\geq 0}$ of \eqref{sfde} solves the following
integral equation
\begin{equation}\label{int}
X(t)=e^{-\lambda_0t}X_0+\int_0^t e^{-\lambda_0(t-s)}Z(X_s)\:\dd
s+J^{\lambda_0}(t),
\end{equation}

\noindent where $J^{\lambda_0}(t)\defeq \int_0^t
e^{-\lambda_0(t-s)}\:\dd W_s$ is the Ornstein-Uhlenbeck process
solving
\begin{equation*}
\left\{
\begin{array}{ll}
\dd u(t)= v(u(t))\dd t + \dd W_t,\quad t\geq 0 \\
u(0)=0.
\end{array}
\right.
\end{equation*}

\noindent Let $\delta>0$. By using \eqref{int} there exists a
positive constant $c_{\delta}>1$ such that for $t\geq 0$,

\begin{equation}\label{x}
|X(t)|^2\leq c_{\delta}
e^{-2\lambda_0t}|X_0|^2+c_{\delta}\frac{M^2}{\lambda_0^2}+(1+\delta)|J^{\lambda_0}(t)|^2
\end{equation}

\noindent On the other hand we know that $\sup\limits_{t\geq
0}\E\Big( e^{\ve|J^{\lambda_0}(t)|^2}\Big)<+\infty$ whenever
$\ve<\lambda_0$. Now Theorem 12.1  in \cite{lifshits} implies that
$$
\sup\limits_{t\geq 0}\E\Big(
e^{\ve\|J_t^{\lambda_0}\|^2}\Big)<+\infty\quad\mbox{for
$\ve<\lambda_0$}.
$$

Therefore, from \eqref{x} and by assuming that the initial condition
$X_0$ is deterministic we have
\begin{equation*}
\sup\limits_{t\geq 0}\E\Big(
e^{\ve\|X_t\|^2}\Big)<+\infty\quad\mbox{for
$\ve(1+\delta)<\lambda_0$}.
\end{equation*}

\noindent From the arbitrariness of $\delta$ we obtain
\begin{equation}\label{moment}
\sup\limits_{t\geq 0}\E\Big(
e^{\ve\|X_t\|^2}\Big)<+\infty\quad\mbox{for $\ve<\lambda_0$}.
\end{equation}

\noindent This implies in particular tightness of the family
$\{X_t:\, t\geq 0\}$ and hence by using a similar argument as in
\cite{Es-vGa-Sch} we deduce the existence of an invariant measure
$\mu$ for the segment process $(X_t)_{t\geq 0}$ on the space $\C$.
Moreover by \cite[Theorem 3]{Sch} we have $\mathcal{L}(X_t)$
converges to $\mu$ in total variation as $t\to +\infty$. Hence
inequality \eqref{moment} yields
$$
\int_{\C}e^{\ve\|x\|^2}\:\mu(\dd x)<+\infty\quad \mbox{for
$\ve<\lambda_0$}.
$$

\noindent Thus, choosing $\lambda$ such that
$4(2L+rL^2)<\lambda<\lambda_0$ yields the integrability condition
\eqref{intcond}. It is classical that the Harnack inequality
\eqref{harnack} implies that the semigroup $(P_t)_{t\geq 0}$ is
strong Feller and irreducible, hence uniqueness of $\mu$ follows
from the classical Doob's Theorem \cite[Theorem 4.2.1]{DaZ96}.
Alternatively, uniqueness follows from \cite{hairer-2009} and
\cite[Theorem 3]{Sch}.

\noindent

\section{Proofs}

\textit{Proof of Theorem \ref{mainthm}}.  As in \cite{MR2215664} we shall employ a coupling argument. Let $x, y \in \C$ be given and  let $X$ denote the solution of $\eqref{sfde}$ starting from initial condition $x \in \C$.  Fix $1>\epsilon >0$ and define  $H : \R^d \to \R^d$,
\[ H(x) = \left\{
\begin{array}{ll}
\frac x {|x|} |x|^\epsilon,  &  \mbox{if } x\ne 0\\
0 &   \mbox{if }x=0.
\end{array}\right.
\]
$H$ is continuous and the gradient of the convex function $ \frac 1 {1+\epsilon} |x|^{1+\epsilon}$ on $\R^d$, hence it is also monotone, i.e.\
\[\langle H(x)-H(y), x-y\rangle \geq 0 \quad \forall x,y \in \R^d,\]
where $\langle \cdot ,\cdot \rangle$ denotes the standard inner product on $\R^d$. Thus, for fixed $\gamma >0$, by the general existence and uniqueness result for monotone SDEs, c.f.\ e.g\ \cite{MR1091217},  there exists a unique process $(\tilde  Y(t))_{t\geq 0}$ solving
 \begin{equation}
\label{csde} \left\{
\begin{aligned}
\dd \tilde Y(t) & = v(\tilde Y(t)) \dd  t + Z(X_t)\,\dd t -\gamma \cdot H(\tilde Y(t)-X(t))\dd t  + \dd W(t),\\
\tilde Y(0)& = y(0),\\
 \end{aligned}
\right.
\end{equation}
and which we  extend by $\tilde Y(t)=y(t)$ for $t \in [-r,0[$.  \\

In particular, in view of \eqref{weakdissip}, for $R(t) = X(t) - \tilde Y(t) $ and $t \geq 0$,
\[  d |R|^2(t)\leq  -2\gamma  \cdot |R|^{1+\epsilon}(t) \dd t ,
\]
i.e.\
\begin{equation} |R(t)|^2   \leq  \bigl(|R(0)|^{{1-\epsilon} } - \gamma(1-\epsilon) \cdot t\bigr)_+^{2/(1-\epsilon)}.
\label{distproc}
\end{equation}

such that  $R(t) =0$ for $t\geq  |R(0)|^{ {1-\epsilon} }/(\gamma (1-\epsilon))$. Hence, for $s \in ]r,T]$, choosing  $\gamma=\gamma_s = |R(0)|^{ {1-\epsilon} }/((s-r)(1-\epsilon))$ implies    $X_t= \tilde Y_t$ in $\C$ for all $t \geq s$.\\

On the other hand we may rewrite equation \eqref{csde} with $\gamma =\gamma_s$ as

 \begin{equation}
\label{gsde}
\left\{
\begin{aligned}
\dd \tilde Y(t)  & = V(\tilde Y_t)\,\dd t    + \dd \tilde W(t), \\
 \tilde Y_0 &=y,
 \end{aligned}
\right.
\end{equation}
where $\dd \tilde W(t) =  \dd W(t) - \zeta(t) \dd t $ with $\zeta (t) = \gamma_s \cdot H(\tilde Y(t)-X(t)) -( Z(X_t)-Z(\tilde Y_t)) $. \\

Now, due to \eqref{distproc} and the Lipschitz bound on $Z$ it holds that
\begin{align}
\frac 1 2 \int_0^T  | \zeta(u)|^2 \dd u  &\leq \gamma_s^2 \int_0^s |R(u)|^{2\epsilon} \dd u +   L^2 \int_0^s \|R_u\|^2 \dd u  \nonumber \\
& \leq   \frac{ |R(0)|^2}{(1-\epsilon^2)(s-r)} +   L^2 \int_0^s \|R_u\|^2 \dd u  \nonumber \\
& \leq  \frac{ |R(0)|^2}{(1-\epsilon^2)(s-r)} +  L^2\, s \, \|R_0\|^2  \quad \mathbb P \mbox{-a.s.},\label{gest}
\end{align}
where in the last step we used the a.s.\ monotonicity of  $u \mapsto |R(u)|$ for $u\geq 0$. \\

In particular, the Novikov condition is satisfied for the exponential martingale $\mathcal E(\xi_.)$ with $\xi(t) = \int_0^t \zeta(s) \dd W(s)$, $t \in [0,T]$, and  by   the Girsanov theorem $(\tilde W_t)_{t\in [0,T]}$  is a Brownian motion under the probability measure $\dd \mathbb Q = D \dd \mathbb P$ with $D= \exp(\int_0^T \langle \zeta(u), \dd W(u)\rangle  - \frac  1 2 \int_0^T |\zeta(u)|^2 \dd u)$, i.e. $(\mathbb Q, \tilde Y)$ is weak solution to \eqref{sfde} starting from $y$. \\

Finally,  for $p>1 $ and $q = p/(p-1)$
\begin{align*}
P_Tf(y) & = \mathbb E_{\mathbb Q} [f(\tilde Y_T)]  =  \mathbb E_{\mathbb P} [D \cdot f(\tilde Y_T)] =   \mathbb E_{\mathbb P} [D \cdot f(  X_T)] \\
& \leq (\mathbb E_{\mathbb P} [D^q])^{\frac 1 q}  (\mathbb E_{\mathbb P} f^p( X_T))^{\frac 1 p} = (\mathbb E_{\mathbb P} [D^q])^{\frac 1 q} \cdot (P_T (f^p)(x))^{\frac 1 p},
\end{align*}
with
\begin{align*}
 \mathbb E_{\mathbb P}( D^q)= \mathbb E_{\mathbb P} \bigl( \exp\bigl(q\int_0^T \zeta(u) dW(u) - \frac q 2 \int_0^T |\zeta(u)|^2\dd u \bigr)\bigr)
 \leq \| \exp({  \frac {q(q-1)}{2} \int_0^T |\zeta(u)|^2 \dd u })\|_{L^{\infty}(\mathbb P)} .
\end{align*}
Hence, due to \eqref{gest}, we arrive at
\[(P_T f(y) )^p \leq P_T(f^p)(x) \exp\Bigl( \frac{ p}{p-1} \Bigl ( \frac  {|x(0)-y(0)|^2} {(s-r)(1-\epsilon^2)} +  s L^2\|x-y\|^2\Bigr)\Bigr) , \]
such that the claim follows by letting $\epsilon \to 0$ and optimizing over $s\in ]r, T]$. \hfill $\Box$ \\

\textit{Proof of Corollary  \ref{stronfell}.}  For  $T>r$  and $x, y \in \C$, proceed as in the previous proof by choosing $\epsilon >0$ and $s=T$. Then for   $f: \C\to \R$ bounded measurable
\begin{align}
|P_Tf(x) -P_Tf(y)| &= | \mathbb  E_{\mathbb Q} f(\tilde Y_T)- \mathbb E _{\mathbb P} f(X_T)| =  |\mathbb E_{\mathbb P} ((1-D) f(X_T))| \nonumber \\
& \leq \|f\|_\infty \sqrt {\mathbb E_{\mathbb P}(|1-D|^2)}= \|f\|_\infty \sqrt{ \mathbb E_{\mathbb P} (D^2)-1}\nonumber \\
& \leq  \|f\|_\infty \sqrt{ \exp \bigl( 2   \bigl( \frac  {|x(0)-y(0)|} {(T-r)(1-\epsilon^2)} +  T \, L^2\|x-y\|^2\bigr)\bigr)-1},  \nonumber
\end{align}
which tends to zero, even uniformly, for $x \to y $ in $\C$. \hfill $\Box$

%

\end{document}